\documentclass[12pt,twoside,final,psamsfonts]{amsart}

\usepackage[psamsfonts]{amssymb}
\usepackage{times,a4wide}

\usepackage{color}		
\usepackage{epsfig}

\theoremstyle{plain}
\newtheorem*{theorem*}{Theorem}

\newtheorem{theorem}{Theorem}[section]

\newtheorem*{proposition*}{Proposition}

\newtheorem*{corollary*}{Corollary}

\newtheorem*{lemma*}{Lemma}

\theoremstyle{definition}
\newtheorem{remark}[theorem]{Remark}
\newtheorem*{remark*}{Remark}

\theoremstyle{definition}
\newtheorem{definition}[theorem]{Definition}
\newtheorem*{definition*}{Definition}
\newcommand{\nc}{\newcommand}

\newcommand{\N}{{\mathbb N}}
\newcommand{\DD}{{\mathbb D}}

\newcommand{\T}{{\mathbb{T}}}

\newcommand{\Int}{\operatorname{Int}}

\nc{\Intl}{\Int_{l^{\infty}}}

\newcommand{\Llra}{\Longleftrightarrow}

\newcommand{\vp}{\varphi}
\nc{\bea}{\begin{eqnarray}}
\nc{\eea}{\end{eqnarray}}
\nc{\beqa}{\begin{eqnarray*}}
\nc{\eeqa}{\end{eqnarray*}}
\nc{\Hi}{H^{\infty}}
\nc{\loi}{\ell^{\infty}}
\nc{\NL}{N^+\vert \Lambda}
\nc{\liL}{\lambda\in\Lambda}
\nc{\nn}{\nonumber}
\nc{\hf}{{\mathcal H}_{\Phi}}
\nc{\hF}{{\mathcal H}_{\Phi}}

\newenvironment{proof*}{\vskip 2mm\noindent {}}{$\blacksquare$ \vskip 2mm}
\numberwithin{equation}{section}

\nc{\card}{\operatorname{card}}
\nc{\tsn}{\tilde{\sigma}_n}
\nc{\tsk}{\tilde{\sigma}_k}
\nc{\tskp}{\tilde{\sigma}_k^+}
\nc{\tskm}{\tilde{\sigma}_k^-}
\nc{\dst}{\displaystyle}
\nc{\vt}{\vartheta}

\title{Extremely weak interpolation in $\Hi$}

\author{Andreas Hartmann}

\address{Equipe d'Analyse,
Institut de Math\'ematiques de Bordeaux
Universit\'e Bordeaux I, 351 cours de la Lib\'eration,
33405 Talence, France}
\email{hartmann@math.u-bordeaux.fr}

\date{\today}

\thanks{This project was elaborated while the author was Gaines Visiting Chair 
at the University of Richmond, also partially supported by the french ANR-project FRAB}

\keywords{Hardy spaces, interpolating sequences, weak interpolation}

\subjclass{30E05, 32A35}

\begin{document}

\begin{abstract}
Given a sequence of points in the unit disk,
a well known result due to Carleson states that if given any point of the sequence
it is possible to
interpolate the value one in that point and zero in all the other points of the sequence,
with uniform control of the norm in the Hardy space of bounded analytic
functions on the disk, then the sequence is an interpolating sequence
(i.e.\ every bounded sequence of values can be interpolated by functions in the
Hardy space). 
It turns out that such a result holds in other spaces.
In this short
note we would like to show that for a given sequence it is sufficient to 
find just {\bf one} function interpolating suitably zeros and ones to deduce
interpolation in the Hardy space.
\end{abstract}

\maketitle

\section{Introduction}

The Hardy space $\Hi$ of bounded analytic functions 
on $\DD$ is equipped with the usual norm $\|f\|_{\infty}=\sup_{z\in\DD}|f(z)|$.
A sequence $\Lambda=\{\lambda_n\}_n\subset \DD$ of points in the unit disk
is called interpolating for $\Hi$, noted $\Lambda\in \Int \Hi$, 
if every bounded sequence of values $v=(v_n)_n\in l^{\infty}$
can be interpolated by a function in $\Hi$. Clearly, for $f\in \Hi$, the sequence $(f(\lambda_n))_n$
is bounded. Hence
\[
 \Lambda\in\Int\Hi \quad \stackrel{\mbox{def}}{\Llra} \quad \Hi|\Lambda=l^{\infty}
\]
(we identify the trace space with a sequence space).
The sequence $\Lambda$ is said to satisfy the Blaschke condition if
$\sum_n (1-|\lambda_n|)<\infty$. In that case, the Blaschke product
$B=\prod_n b_{\lambda_n}$, where $b_{\lambda}(z)=\frac{|\lambda|}{\lambda}
\frac{z-\lambda}{1-\overline{\lambda}z}$ is the normalized M\"obius transform
($\lambda\in\DD$),
converges uniformly on every compact set of $\DD$ to a function in $\Hi$ with
boundardy values
$|B|=1$ a.e.\ on $\T$. Carleson proved  (see \cite{carl}) that 
\[
 \Lambda\in \Int \Hi\quad \Llra \quad \inf_n |B_n(\lambda_n)|=\delta>0,
\]
where $B_n=\prod_{k\neq n} b_{\lambda_k}$. The latter condition will be termed
Carleson condition, and we shall write $\Lambda\in (C)$ when $\Lambda$
satisfies this condition. Carleson's result can be reformulated 
using the notion of weak interpolation.

\begin{definition}
A sequence $\Lambda$ of points in $\DD$ is called a weak interpolating sequence in $\Hi$,
noted $\Lambda\in \Int_w\Hi$,
if for every $n\in \N$ there exists a function $\vp_n\in \Hi$ such that
\begin{itemize}
\item for every $n\in \N$, $\vp_n(\lambda_k)=\delta_{nk}$,
\item $\sup_n\|\vp_n\|_{\infty}<\infty$.
\end{itemize}
\end{definition}

Now, when $\Lambda\in (C)$, setting $\vp_n=B_n/B_n(\lambda_n)$ we obtain a family
of functions satisfying the conditions of the definition. Hence $\Lambda\in \Int_w\Hi$.
And from Carleson's theorem we get 
\[
 \Lambda\in \Int\Hi\quad \Llra \quad \Lambda\in \Int_w\Hi.
\]

With suitable definitions of interpolating and weak interpolating sequences, such a result
has been shown to be true in Hardy spaces $H^p$ (see \cite{ShHSh} for $1\le p<\infty$
and \cite{Kab} for $0<p<1$) as well as in Bergman spaces (see \cite{SchS1}) and in certain
Paley-Wiener and Fock spaces (see \cite{SchS2}).

One also encounters the notion of ``dual boundedness'' for such sequences (see \cite{Am}),
and in a suitable context it is related to so-called uniform minimality of sequences
of reproducing kernels (see e.g.\ \cite[Chapter C3]{Nik02} for some general facts).

Using a theorem by Hoffmann we want to show here that given a separated sequence
$\Lambda$, then there is a splitting of $\Lambda=\Lambda_0\cup\Lambda_1$ such 
that if there exists {\bf just one} function $f\in \Hi$ vanishing on $\Lambda_0$ and being 1 on
$\Lambda_1$, then the sequence is interpolating in $\Hi$.

The author does not claim that such a result is anyhow useful to test whether a sequence is interpolating
or not,
but that it might be of some theoretical interest.

\section{The result}

Let us begin by recalling Hoffman's result (which can e.g.\ be found in 
Garnett's book, \cite{Gar}):
\begin{theorem}[Hoffman]
For $0<\delta<1$, there are constants $a=a(\delta)$ and $b=b(\delta)$
such that the Blaschke product $B(z)$ with zero set
$\Lambda$ has a nontrivial factorization
$B=B_0B_1$ and
\beqa
 a|B_0(z)|^{1/b}\le |B_1(z)|\le \frac{1}{a}|B_0(z)|^b
\eeqa
for every $z\in \DD\setminus \bigcup_{\lambda\in\Lambda}D(\lambda,\delta)$,
where $D(\lambda,\delta)=\{z\in\DD:|b_{\lambda}(z)|<\delta\}$ 
is the pseudohyperbolic disk.
\end{theorem}

In view of this theorem, 
given any Blaschke sequence of points $\Lambda$ in the disk and
a constant $\delta\in (0,1)$, we will set $\Lambda_0$
to be the zero set of $B_0$ and $\Lambda_1$ to be the zero set of $B_1$ where $B=B_0B_1$
is a Hoffman factorization of $B$. We will refer to $\Lambda=\Lambda_0\cup\Lambda_1$
as a {\it $\delta$-Hoffman decomposition} of $\Lambda$. 

Recall that a sequence $\Lambda$ is separated if there exists a constant $\delta_0>0$ such that
for every $\lambda,\mu\in \Lambda$, $\lambda\neq \mu$, $|b_{\lambda}(\mu)|\ge \delta_0$.
For such a sequence, we will call
a {\it corresponding} Hoffman decomposition a Hoffman decomposition
associated with $\delta=\delta_0/2$.

\begin{theorem}\label{thm2.2}
A separated sequence $\Lambda$ in the unit disk with corresponding Hoffman decomposition
$\Lambda=\Lambda_0\cup\Lambda_1$ is interpolating for $\Hi$ if and only
if there exists a function $f\in \Hi$ such that $f|\Lambda_0=0$ and 
$f|\Lambda_1=1$.
\end{theorem}

The condition is clearly necessary.

\begin{proof}[Proof of Theorem]
Preliminary observation: by 
factorization in $\Hi$ (see e.g.\ \cite{Gar}), we have $f=B_0F$ where $F$ is a bounded
analytic functions (that could contain inner factors). Then
for every $\mu\in \Lambda_1$
\[
 1=f(\mu)=|B_0(\mu)||F(\mu)|\le c|B_0(\mu)|
\]
which shows that 
\[
 |B_0(\mu)|\ge \eta:=1/c.
\]
Replacing $f$ by $g=1-f$ we obtain a function vanishing now on $\Lambda_1$ and
being 1 on $\Lambda_0$. And the same argument as before shows that
for $\mu\in\Lambda_0$
\[
 |B_1(\mu)|\ge \eta
\]
(let us agree to use the same $\eta$ here).

Pick now $\mu\in \Lambda_1$. Then 
\[
 |B_0(\mu)|\ge \eta.
\]
We have to check whether such an estimate holds also for the second piece.
Now, let $z\in \partial D(\mu,\delta)$ (note that $\delta=\delta_0/2$, where $\delta_0$ is the
separation constant of $\Lambda$, so that this disk
is far from the other points of $\Lambda$). Then by Hoffman's theorem
\[
 |B_{\Lambda_1}(z)|\ge a|B_0(z)|^{1/b}
\]
Hence
\[
 |B_{\Lambda_1\setminus\{\mu\}}(z)||b_{\mu}(z)|\ge a|B_0(z)|^{1/b}
\]
and
\[
 |B_{\Lambda_1\setminus\{\mu\}}(z)|\ge \frac{a}{\delta}|B_0(z)|^{1/b}
\]
Now $B_{\Lambda_1\setminus\{\mu\}}$ and $B_0$ do not vanish in $D(\mu,\delta)$.
We thus can take powers of $B_0$ and divide through getting a function
$B_{\Lambda_1\setminus\{\mu\}}/B_0^{1/b}$ not vanishing in $D(\mu,\delta)$. By the
minimum modulus principle we obtain
\[
 \left| \frac{B_{\Lambda_1\setminus\{\mu\}}(z)}{B_0^{1/b}(z)}\right|\ge \frac{a}{\delta}
\]
for every $z\in D(\mu,\delta)$ and especially in $z=\mu$ so that
\[
 |B_{\Lambda_1\setminus\{\mu\}}(\mu)|\ge \frac{a}{\delta}\eta^{1/b}.
\]
Hence
\[
 |B_{\Lambda\setminus\{\mu\}}(\mu)|=|B_0(\mu)||B_{\Lambda_1\setminus\{\mu\}}(\mu)|
 \ge \frac{a}{\delta}\eta^{1+1/b}.
\]

By the preliminary observation above,
the same argument can be carried through when $\mu\in\Lambda_0$, so that 
for every $\mu\in\Lambda$ we get
\[
 |B_{\Lambda\setminus\{\mu\}}(\mu)|\ge c
\]
for some suitable $c>0$. Hence $\Lambda\in(C)$ and we are done.
\end{proof}

\begin{remark}
1) It is clear from the proof that it is sufficient that
there is an $\eta>0$ with
\bea\label{eq1}
 \inf_{\mu\in \Lambda_1}|B_0(\mu)|\ge\eta\quad \mbox{and}\quad
 \inf_{\mu\in \Lambda_0}|B_1(\mu)|\ge\eta
\eea
This means that in terms of Blaschke products, we need {\bf two} functions instead
of the sole function $f$ (which is of course not unique) as stated in the theorem.
One could raise the question whether it would be sufficient to have only one of
the conditions in \eqref{eq1} (the condition is clearly necessary).
Suppose we had the first condition
\[
 \inf_{\mu\in \Lambda_1}|B_0(\mu)|\ge\eta
\]
Then in order to obtain the condition of the theorem, we would need to 
multiply $B_0$ by a function $F\in\Hi$ such that $(B_0F)(\mu)=1$ for
every $\mu\in \Lambda_1$. In other words we need that $B_0+B_1\Hi$ is
invertible in the quotient algebra $\Hi/B_1\Hi$ under the condition 
that $0<\eta\le |B_0(\mu)|\le 1$. This is possible when $\Lambda_1$
is a finite union of interpolating sequences in $\Hi$ (which in our case
boils down to interpolating sequences since we have somewhere assumed
that $\Lambda$, and hence $\Lambda_1$, is separated). See for example
\cite{AHth} for this, but it can also be deduced from Vasyunin's earlier
characterization of the trace of $\Hi$ on finite union of interpolating
sequences (see \cite{vas}).

We do not know the general answer to this invertibility problem
when $\Lambda_1$ is not assumed to be a finite union of interpolating sequences. 

2) Another question that could be raised is whether in Theorem \ref{thm2.2} the
assumption of being separated can be abandoned. At least Hoffman's theorem
does not allow us to deduce that the sequence is separated. As an example,
one could have a union of two interpolating sequence the elements of which come
arbitrarily close to each other. Write $\Lambda=\bigcup_n \sigma_n$ where $\sigma_n$
contains two close points of $\Lambda$ one of which is of the first interpolating 
sequence and the other one from the second interpolating sequence. Let
$\Lambda_0$ be the union of the even indexed $\sigma_n$'s and $\Lambda_1$
the odd indexed $\sigma_n$'s we obtain a Hoffman decomposition for which
we can find $f$ as in the theorem, but $\Lambda$ is not interpolating. 
\end{remark}

\end{document}